\documentclass[11pt]{article}

\usepackage{amsfonts}
\usepackage{amsmath}

\newcommand{\beq}{\begin{equation}}
\newcommand{\eeq}{\end{equation}}
\newcommand{\bea}{\begin{eqnarray}}
\newcommand{\eea}{\end{eqnarray}}
\newcommand{\beas}{\begin{eqnarray*}}
\newcommand{\eeas}{\end{eqnarray*}}

\newtheorem{theorem}{Theorem}[section]

\newtheorem{proposition}[theorem]{Proposition}
\newtheorem{corollary}[theorem]{Corollary}
\newtheorem{lemma}[theorem]{Lemma}
\newtheorem{remark}[theorem]{Remark}
\newtheorem{example}[theorem]{Example}
\newtheorem{examples}[theorem]{Examples}
\newtheorem{foo}[theorem]{Remarks}

\newenvironment{proof}{\addvspace{\medskipamount}\par\noindent{\it
Proof}.}
{\unskip\nobreak\hfill$\Box$\par\addvspace{\medskipamount}}

\newcommand{\p}{\partial}

\newcommand{\bM}{\mathbb M}

\newcommand{\M}{\mathbb M}

\newcommand{\R}{\mathbb R}
\newcommand{\f}{e^{\lambda f}}
\newcommand{\pl}{e^{\lambda \psi}}

\title{Perelman's entropy and doubling property on Riemannian manifolds}

\author{Fabrice Baudoin, Nicola Garofalo}

\begin{document}

\maketitle

\begin{abstract}
The purpose of this work is to study some monotone functionals of the heat kernel on a complete Riemannian manifold with nonnegative Ricci curvature. In particular, we show that on these manifolds, the gradient estimate of Li and Yau \cite{LY}, the gradient estimate of Ni \cite{Ni}, the monotonicity of the Perelman's entropy and the volume doubling property are all consequences of an entropy inequality recently discovered in \cite{baudoin-garofalo}. The latter is a linearized version of a logarithmic Sobolev inequality that is due to D. Bakry and M. Ledoux \cite{Bakry-Ledoux}. 
\end{abstract}

\section{Introduction}\label{S:intro}

Comparison theorems play a fundamental role in Riemannian geometry. Two results which are particularly important are the Bishop-Gromov volume comparison theorem and the Laplacian comparison theorem. We will focus on manifolds with nonnegative Ricci tensor. In such context the volume comparison theorem states that, given a complete $n$-dimensional Riemannian manifold, then  for every $x\in \M$ the quotient $\frac{\text{Vol}_n(B(x,r)}{r^n}$ 
is a non-increasing function of $r>0$. In particular, this implies the following global doubling condition
\begin{equation}\label{d}
\text{Vol}_n(B(x,2r))\le 2^n \text{Vol}_n(B(x,r)),\ \ \ \ x\in \M,\ r>0.
\end{equation}
The Laplacian comparison theorem  (which is a special case of the Hessian comparison theorem) states instead that, given a fixed base point $x_0\in \M$, at any point $x\in \M$ where the Riemannian distance $d(x) = d(x,x_0)$ is smooth one has
\begin{equation}\label{lap}
\Delta d(x) \le \frac{n-1}{d(x)},
\end{equation}
where $\Delta$ indicates the Laplace-Beltrami operator on $\M$. Furthermore, such inequality holds globally on $\M$ in the sense of distributions.
To cite from \cite{CLN}: ``The ideas behind these elementary results have a profound influence on geometric analysis and Ricci flow". The inequalities \eqref{d} and \eqref{lap} are also intimately connected to the celebrated inequality of Li and Yau \cite{LY} which, for a complete manifold with Ricci $\ge 0$, states that for every smooth nonnegative function $f$ on $\M$ one has 
\begin{equation}\label{LY}
\Delta(\ln P_t   f) \ge - \frac{n}{2t},
\end{equation}
where $P_t = e^{t\Delta}$ indicates the heat semigroup associated with $\Delta$. The Li-Yau inequality, in turn, is deeply related to an entropy inequality for the  Ricci flow, more recently discovered by Perelman, see \cite{Per}. 

The purpose of this work is to show that the above cited results may actually be seen as consequences of a new linearized version, found in \cite{baudoin-garofalo}, of a logarithmic Sobolev inequality that is due to Bakry and Ledoux \cite{Bakry-Ledoux}. In particular, we will obtain a new proof of the volume doubling property (\ref{d}) that does not rely on the Laplacian comparison theorem (\ref{lap}). It is our hope, and the primary motivation of the present work, that the techniques developed here may apply to situations, such as complete sub-Riemannian manifolds, where a Laplacian comparison theorem is not available. On the other hand, it has been recently shown in \cite{baudoin-garofalo}, that on such manifolds there is a notion of Ricci tensor which generalizes its Riemannian predecessor, and that one can successfully develop a theory that parallels the Riemannian one. 

Finally, we observe that as a by-product of our work we obtain seemingly new global heat kernel bounds on complete Riemannian manifolds with strictly positive Ricci curvature.

\section{Li-Yau type estimates and Perelman's entropy}

Let $\mathbb{M}$ be a complete Riemannian manifold of dimension $n$. We denote by $\Delta$ the Laplace-Beltrami operator on $\mathbb{M}$. Throughout the paper we use the convention that $\Delta$ is a nonpositive operator.
Let us recall that the Riemannian measure is symmetric and invariant for $\Delta$ defined on the space $ C^\infty_0(\bM)$ of smooth and compactly supported functions. More precisely,  for $f,g \in C^\infty_0(\bM)$, one has
\[
\int_{\bM} f \Delta g d\mu=-\int_{\bM} (\nabla f, \nabla g) d\mu=\int_{\bM} g \Delta f d\mu,
\]
and 
\[
\int_{\bM} \Delta f d\mu=0.
\]
The assumption that the metric on $\bM$ is complete implies that $\Delta$ is essentially self-adjoint on  $ C^\infty_0(\bM)$ (see \cite{davies}). The Friedrichs extension of $\Delta$ is then still denoted by $\Delta$.

The operator $\Delta$ is the generator of a strongly continuous contraction  semigroup $(P_t)_{t \ge 0}$ on $L^p (\bM,\mu)$, $1 \le p \le +\infty$, that moreover satisfies the sub-Markov property: If $u \in L^2 (\bM,\mu)$ is such that $0 \le u \le 1$, then $0 \le P_t u \le 1$. Actually, due to the ellipticity of $\Delta$, $(t,x) \rightarrow P_t f(x)$ is
smooth on $\mathbb{M}\times (0,\infty) $ and
\[ P_t f(x)  = \int_{\mathbb M} p(x,y,t) f(y) d\mu(y),\ \ \ f\in
C^\infty_0(\mathbb M),\] where $p(x,y,t) > 0$ is the so-called heat
kernel associated to $P_t$. Such function is smooth outside the
diagonal of $\bM\times \bM$, and it is symmetric, i.e., \[ p(x,y,t)
= p(y,x,t). \]
 By the
semi-group property for every $x,y\in \bM$ and $0<s,t$ we have
\begin{align*}
p(x,y,t+s) & = \int_\bM p(x,z,t) p(z,y,s) d\mu(z) = \int_\bM p(x,z,t)
p(y,z,s) d\mu(z)
\\
& = P_s(p(x,\cdot,t))(y).
\end{align*}

\subsection{The curvature dimension inequality}\label{SS:CD}

Associated to $\Delta$ are the two following differential bilinear forms defined on the space of smooth functions:
\begin{equation*}
\Gamma(f,g) =\frac{1}{2}\big[\Delta(fg)-f\Delta g-g\Delta f\big]=( \nabla f
, \nabla g),
\end{equation*}
and
\begin{equation*}
\Gamma_{2}(f,g) = \frac{1}{2}\big[\Delta\Gamma(f,g) - \Gamma(f,
\Delta g)-\Gamma (g,\Delta f)\big],
\end{equation*}
where $f, g :\mathbb{M} \rightarrow \mathbb{R}$ are smooth
functions. Henceforth, we will simply write $\Gamma(f) = \Gamma(f,f)$, $\Gamma_2(f) = \Gamma_2(f,f)$. As an application of the Bochner's formula which we
can re-write as \[ \Delta \Gamma(f,f) = 2 ||\nabla^2 f||^2 + 2
\Gamma(f,\Delta f) + 2\ \text{Ric}(\nabla f,\nabla f), \] one
obtains
\[
\Gamma_{2}(f,f)= \| \nabla^2 f \|_2^2 +\text{Ric}(\nabla f, \nabla
f).
\]
With the aid of Schwarz inequality, which gives $\| \nabla^2 f
\|_2^2\ge \frac{1}{n} (\Delta f)^2$, the assumption that the Riemannian Ricci tensor on $\bM$ is bounded from below by $\rho \in \mathbb{R}$,  translates then into the so-called curvature-dimension inequality $\mathbf{CD}(\rho,n)$:
\begin{equation}\label{curvature_dimension}
\Gamma_{2}(f,f) \ge \frac{1}{n} (\Delta f)^2 + \rho 
\Gamma(f,f).
\end{equation}

In the remainder of the section, we assume that the inequality $\mathbf{CD}(\rho,n)$ is satisfied, which means that Riemannian Ricci tensor of $\bM$ is bounded from below by $\rho$.

\subsection{The main variational inequality}

In a first step, we translate some of the results obtained by Baudoin-Garofalo \cite{baudoin-garofalo} to this Riemannian framework.
Let us consider a fixed $T>0$. 
For $f\in C^\infty_0(\bM)$ and $f\geq 0$ we introduce the  functional
\[
\Phi (t)=P_t \left( (P_{T-t} f) \Gamma (\ln P_{T-t}f) \right),\ \ \ 0\le t < T.
\]
We first have the following key lemma that is due to Bakry-Ledoux, see \cite{Bakry-Ledoux}.

\begin{lemma}\label{L:derivatives}
For every $0\le t < T$. We have
\[
\int_0^t \Phi (s) ds=P_t ( P_{T-t} f  \ln P_{T-t} f)-P_T f \ln P_Tf
\]
\[
\Phi' (t)=2P_t \left( (P_{T-t} f) \Gamma_2 (\ln P_{T-t}f) \right).
\]
\end{lemma}

We may then use the results of \cite{baudoin-garofalo} to obtain the following fundamental variational inequality, whose proof is reproduced here for the sake of completeness.

\begin{theorem}\label{variational}
Let $a:[0,T] \rightarrow \mathbb R^+$ be a smooth and positive function. We have for every smooth $\gamma :[0,T] \rightarrow \mathbb{R}$,
\begin{align}\label{variational inequality}
(a \Phi )'  & \ge \left(a' -\frac{4a\gamma}{n} +2 \rho a \right)\Phi  +\frac{4a\gamma}{n} \Delta P_{T} f - \frac{2a\gamma^2}{n} P_T f.
\end{align}
\end{theorem}

\begin{proof}
To prove this result we apply the
curvature-dimension inequality $\mathbf{CD}(\rho,n)$, in combination with Lemma
\ref{L:derivatives}. If $a$ is a positive function we 
obtain
\begin{align*}
(a \Phi )'  & \ge  a' \Phi+2a \rho \Phi+\frac{2a}{n} \left( P_t ( (P_{T-t} f) (\Delta \ln P_{T-t} f )^2 )\right)
\end{align*}
Now, for every $\gamma \in \mathbb{R}$ one has
\[
(\Delta \ln P_{T-t} f )^2 \ge 2\gamma \Delta \ln P_{T-t}f -\gamma^2.
\]
Furthermore,
\[
 \Delta \ln P_{T-t}f=\frac{\Delta P_{T-t}f}{P_{T-t}f} -\Gamma(\ln P_{T-t} f ).
 \]
Therefore,
 \begin{align*}
(a \Phi )'  & \ge \left(a' -\frac{4a\gamma}{n}+2\rho a \right)\Phi  +\frac{4a\gamma}{n} \Delta P_{T} f - \frac{2a\gamma^2}{n} P_T f,
\end{align*}
which is the desired conclusion.
\end{proof}

We shall often use Theorem \ref{variational} for the heat kernel itself, that is, we take for $f$ a Dirac mass. If we indicate $p_t(x,y) = p(x,y,t)$, we thus obtain the following result.

\begin{lemma}\label{L:variational_heatkernel}
For any given $T>0$ consider the function \[
\Phi (t)=P_t \left( p_{T-t}  \Gamma (\ln p_{T-t}) \right),
\]
where $0\le t<T$. Let $a:[0,T] \rightarrow \R^+$ be a smooth and positive function. We have for every smooth $\gamma :[0,T] \rightarrow \mathbb{R}$,
\begin{align*}
(a \Phi )'  & \ge \left(a' -\frac{4a\gamma}{n} +2\rho a\right)\Phi  +\frac{4a\gamma}{n} \Delta p_{T}  -\frac{2a\gamma^2}{n} p_T.
\end{align*}

\end{lemma}

\subsection{A family of Li-Yau inequalities}

As a first application of Theorem \ref{variational} we use it to derive a family of Li-Yau type inequalities. We choose the function $\gamma$ in a such a way that
\[
a' -\frac{4a\gamma}{n} +2\rho a=0.
\]
That is
\[
\gamma=\frac{n}{4} \left( \frac{a'}{a}+2\rho \right).
\]
Integrating  the inequality (\ref{variational inequality}) from $0$ to $T$, and denoting $V=\sqrt{a}$, we obtain the following result.

\begin{proposition}\label{Generalized Li-Yau}
Let $V:[0,T]\rightarrow \mathbb{R}^+$ be a smooth function such that
\[
V(0)=1, V(T)=0.
\]
We have
\begin{align}\label{family Li Yau}
\Gamma (\ln P_T f) & \le \left( 1-2\rho\int_0^T V^2(s) ds\right) \frac{\Delta P_Tf}{P_T f}
\\
& +\frac{n}{2} \left(  \int_0^T V'(s)^2 ds +\rho^2  \int_0^T V(s)^2 ds -\rho \right).
\notag\end{align}
\end{proposition}

A first family of interesting inequalities may be obtained with the choice
\[
V(t)=\left( 1-\frac{t}{T}\right)^\alpha, \alpha>\frac{1}{2}.
\]
In this case we have
\[
 \int_0^T V(s)^2 ds=\frac{T}{2\alpha+1}
 \]
 and
\[
 \int_0^T V'(s)^2 ds=\frac{\alpha^2}{(2\alpha-1)T},
 \]
 so that, according to (\ref{family Li Yau}),
 \begin{equation}\label{Li Yau lineaire}
 \Gamma (\ln P_T f) \le \left( 1-\frac{2\rho T}{2\alpha+1}\right) \frac{\Delta P_Tf}{P_T f} +\frac{n}{2} \left(   \frac{\alpha^2}{(2\alpha-1)T}+\frac{\rho^2 T}{2\alpha+1} -\rho \right).
 \end{equation}
In the case, $\rho=0$ and $\alpha=1$,  \eqref{Li Yau lineaire} reduces to the celebrated Li-Yau inequality:
\begin{align}\label{LiYau}
\Gamma(\ln P_t f) \le \frac{  \Delta P_t f }{P_t f }  
+ \frac{n}{2t},\ \ \ \ \ \ t>0 .
\end{align}

Although in the sequel of this paper we will focus on the case $\rho=0$, let us presently discuss the case $\rho >0$. In that case, from Bonnet-Myers theorem the manifold $\mathbb{M}$ has to be compact and we will assume in the following that $\mu (\bM)=1$.

Using the inequality (\ref{Li Yau lineaire}) with $\alpha=3/2$ leads to the Bakry-Qian inequality (see \cite{BQ}):  
\[
 \frac{\Delta P_tf}{P_t f}\le \frac{n \rho}{4},\ \ \ \ \ t \ge \frac{2}{\rho}.
 \]
 Also, by using the inequality (\ref{family Li Yau}) with
 \[
 V(t)=\frac{e^{-\frac{\rho t}{3}} (e^{-\frac{2\rho t}{3}}-e^{-\frac{2\rho T}{3}})}{1-e^{-\frac{2\rho T}{3}}},
 \]
 we obtain the following inequality that does not seem to have explicitly been noted in the literature: 
 \[
\Gamma(\ln P_t f) \le e^{-\frac{2\rho t}{3}}  \frac{  \Delta P_t f }{P_t f } +\frac{n\rho}{3} \frac{e^{-\frac{4\rho t}{3}}}{ 1-e^{-\frac{2\rho t}{3}}},\ \ \ \ \ t \ge 0.
 \]
 Integrating the latter inequality along geodesics gives the following parabolic Harnack inequality for the heat kernel: 
 \[
P_t f(y) \ge P_s f(x) \left( \frac{1-e^{-\frac{2\rho s}{3}}}{1-e^{-\frac{2\rho t}{3}} }\right)^{n/2} e^{-\frac{\rho}{6} \frac{d(x,y)^2}{ e^{\frac{2\rho t}{3}}-e^{\frac{2\rho s}{3} }}  },\ \ \ \ 0\le s<t.
\]
If we now apply the latter inequality to a delta function, and taking into account the facts that $\lim_{t \to +\infty} p_t (x,y)=1$, $\lim_{t \to 0} t^{n/2} p_t(x,x)=\frac{1}{(4\pi)^{n/2}}$, we end up with the following elegant and global bounds for the heat kernel $p_t(x,y)$:
\[
\left( \frac{\rho}{6\pi} \right)^{n/2}  \frac{1}{(1-e^{-\frac{2\rho t}{3}} )^{n/2} }e^{-\frac{\rho}{6} \frac{d(x,y)^2}{ e^{\frac{2\rho t}{3}}-e^{\frac{2\rho s}{3} }}  } \le p_t(x,y) \le  \frac{1}{(1-e^{-\frac{2\rho t}{3}} )^{n/2} }.
\]
Observe that these seemingly new heat kernel bounds contain the geometric bound:
\[
\rho \le 6\pi,
\]
which, recalling our normalization $\mu (\bM)=1$, implies in the general case:
\[
\mu(\bM) \le \left( \frac{6\pi}{\rho} \right)^{n/2}.
\]

\subsection{Monotonicity of the entropy}\label{SS:mono}

Hereafter, we assume that $\rho=0$. That means that we work on a complete Riemannian manifold with nonnegative Ricci curvature. For $0\le t<T$ we now introduce the functional
\[
W_t \overset{def}{=}  t \Gamma(\ln p_t) - 2 t \frac{\Delta p_t}{p_t} -   \ln p_t -
\frac{n}{2} \ln t.
\]

We have the following basic result.

\begin{proposition}\label{P:monotonicity_entropy}
For $0\le t<T$ one has \[
p_T W_T \le P_t (p_{T-t} W_{T-t}).
\]
\end{proposition}
\begin{proof}
To establish    this result we choose
\[
a(t)=T-t,
\]
\[
\gamma(t)=-\frac{n}{2(T-t)},
\]
so that 
\[
a\gamma \equiv - \frac{n}{2},\ \ \ \ a\gamma^2 = \frac{n^2}{4(T-t)}.
\]
With this choice we now integrate the variational inequality of Lemma \ref{L:variational_heatkernel} from 0 to $t <T$, obtaining
\begin{align*}
& (T-t) P_t \left( p_{T-t}  \Gamma \left(\ln p_{T-t}) \right)\right) \ge  T p_T \Gamma\left(\ln p_T\right)  +\int_0^t \Phi(s) ds - 2 t \Delta p_T 
\\
& - \frac{n}{2} p_T \ln T + \frac{n}{2} p_T \ln(T-t)
\\
& = T p_T \Gamma\left(\ln p_T\right) + P_t ( p_{T-t}  \ln p_{T-t})-p_T  \ln p_T - 2T \Delta p_T
\\
& + 2 (T-t) \Delta p_T - \frac{n}{2} p_T \ln T + \frac{n}{2} p_T \ln(T-t).
\end{align*}

We now have
\[
p_T W_T =  T p_T \Gamma(\ln p_T) - 2 T \Delta p_T -  p_T \ln p_T -
\frac{n}{2} p_T \ln T .
\]
Comparing with the previous inequality we obtain
\begin{equation}\label{ine}
(T-t) P_t \left( p_{T-t}  \Gamma \left(\ln p_{T-t}) \right)\right) -  P_t ( p_{T-t}  \ln p_{T-t}) - 2 (T-t) \Delta p_T  - \frac{n}{2} p_T \ln(T-t) \ge p_T W_T.
\end{equation}
Finally, we have
\begin{align*}
P_t(p_{T-t} W_{T-t}) & =  P_t\left( p_{T-t} \left[(T-t) \Gamma(\ln p_{T-t}) - 2 (T-t) \frac{\Delta p_{T-t}}{p_{T-t}} -   \ln p_{T-t} -
\frac{n}{2} \ln (T-t)\right]\right)
\\
& = (T-t) P_t \left( p_{T-t}  \Gamma \left(\ln p_{T-t} \right)\right) - 2 (T-t) \Delta P_t(p_{T-t})
\\
& -  P_t\left(p_{T-t} \ln p_{T-t}\right) - \frac{n}{2}  \ln (T-t) P_t(p_{T-t})
\\
& = (T-t) P_t \left( p_{T-t}  \Gamma \left(\ln p_{T-t} \right)\right) -  P_t\left(p_{T-t} \ln p_{T-t}\right) 
\\
&  - 2 (T-t) \Delta p_T  - \frac{n}{2} p_T \ln (T-t).
\end{align*}
Using this expression we finally have from \eqref{ine}
\[
P_t(p_{T-t} W_{T-t}) \ge p_T W_T, 
\]
thus completing the proof of the proposition.

\end{proof}

The next lemma plays a crucial role. Its proof is based on the asymptotic expansion of Minashisundaram and Pejiel of the heat kernel $p_t(x,y)$.

\begin{lemma}
\[
\lim_{t \to T} \frac{P_t (p_{T-t} W_{T-t})}{p_T} =n + \frac{n}{2} \log(4\pi)
\]
\end{lemma}

As a conclusion, we recover an inequality pointed out by  Lei Ni \cite{Ni}.

\begin{theorem}
For $t>0$,
\[
W_t=  t \Gamma(\ln p_t) - 2 t \frac{\Delta p_t}{p_t} -   \ln p_t -
\frac{n}{2} \ln t\le n + \frac{n}{2} \log(4\pi)
\]
\end{theorem}

So for the Perelman's entropy,
\[
\mathcal{W}_t =\int_{\bM} p_t W_t d\mu=t \int_{\bM} p_t  \Gamma(\ln p_t) d\mu -\int_{\bM} p_t \ln p_t d\mu-\frac{n}{2} \ln t,
\]
we obtain
\[
\mathcal{W}_t \le n + \frac{n}{2} \log(4\pi).
\]

\section{Exponential integrability}

In this section we establish the following crucial result.

\begin{theorem}\label{T:doubling0}
Let $\M$ be a complete $n$-dimensional Riemannian manifold such that Ric$\ge 0$. There exist an absolute constant $K>0$, and $A>0$, depending only on $n$, such that
\begin{equation}\label{doubling}
 P_{Ar^2}(\mathbf 1_{B(x,r)})(x) \ge K, \ \ \ \ \ x\in \M, r>0.
\end{equation}
\end{theorem}

\begin{proof}
We use Theorem \ref{variational} in which we choose
\[
a(t)=\tau+T-t,
\]
\[
\gamma(t)=-\frac{n}{4(\tau+T-t)} 
\]
where  $\tau >0$ will later be optimized. Noting that we presently have
\[
a' \equiv 1,\ \ \ a\gamma \equiv - \frac{n}{4},\ \ \ \ a\gamma^2 = \frac{n^2}{16(\tau + T - t)^2},
\]
we obtain the inequality
\begin{equation}\label{ine2}
\tau P_T(f \Gamma(\ln f)) -(T+\tau) P_T f \Gamma(\ln P_T f)
\ge -T \Delta P_T f -\frac{n}{8} \ln ( 1+\frac{T}{\tau}) P_T f
\end{equation}
In what follows we consider a bounded function $f$ on $\M$ such that 
$\Gamma(f) \le 1$ almost everywhere on $\M$.  For any $\lambda \in \mathbb R$
we consider the function $\psi$ defined by 
\[
\psi(\lambda,t) = \frac{1}{\lambda} \log P_t(e^{\lambda f}), \ \ \ \text{or alternatively}\ \ \ P_t(e^{\lambda f}) = e^{\lambda \psi}.
\]
Notice that Jensen's inequality gives
\[
\lambda \psi = \log(\pl) = \log P_t(\f) \ge P_t(\log \f) = \lambda P_t f,
\]
and so we have
\begin{equation}\label{Jensen}
P_t f \le \psi.
\end{equation}

We now apply \eqref{ine2} to the function $e^{\lambda f}$, obtaining
\begin{align*}
& \lambda^2 \tau P_T\left(\f \Gamma(f)\right) - \lambda^2 (T+\tau) \pl \Gamma(\psi)
\\
& \ge - T \Delta P_T(\f) - \frac{n}{8} \pl \ln\left(1+\frac{T}{\tau}\right).
\end{align*}
Keeping in mind that 
$\Gamma(f) \le 1$, we see that 
\[
P_T(\f \Gamma(f)) \le P_T(\f) = \pl.
\]
Using this observation in combination with the fact that
\[
\Delta \left(P_t (e^{\lambda f})\right) = \frac{\p}{\p t} \left(P_t (e^{\lambda f})\right) = \frac{\p \pl}{\p t} = \lambda\pl \frac{\p \psi}{\p t} ,
\]
and switching notation from $T$ to $t$, we infer
\begin{align*}
& \lambda^2 \tau  \ge \lambda^2 (t+\tau) \pl \Gamma(\psi)- \lambda t \frac{\p \psi}{\p t} - \frac{n}{8}  \ln\left(1+\frac{t}{\tau}\right).
\end{align*}
The latter inequality finally gives
\begin{equation}\label{ine3}
\frac{\p \psi}{\p t} \ge - \frac{\lambda}{t}\left(\tau + \frac{n}{8\lambda^2} \ln\left(1+\frac{t}{\tau}\right)\right).
\end{equation}
We now optimize the right-hand side of \eqref{ine3} with respect to $\tau$. We notice explicitly that the maximum value of the right-hand side is attained at
\[
\tau_0 = \frac{t}{2} \left(\sqrt{1+\frac{n}{2\lambda^2t}} - 1\right).
\]
If we substitute  such value in \eqref{ine3} we find
\begin{equation}\label{ine4}
-\frac{\p \psi}{\p t} \le \frac{\lambda}{2}\left(\sqrt{1+\frac{n}{2\lambda^2 t}} - 1\right) + \frac{n}{8\lambda t} \ln\left(1+\frac{2}{\sqrt{1+\frac{n}{2\lambda^2 t}} - 1}\right) = \lambda G\left(\frac{1}{\lambda^2 t}\right),
\end{equation}
where we have set
\[
G(s) = \frac{1}{2}\left(\sqrt{1+\frac{n}{2}s} - 1\right) + \frac{n}{8} s \ln\left(1 + \frac{2}{\sqrt{1+\frac{n}{2} s} - 1}\right),\ \ \ \ s>0.\]
Notice that $G(s) \to 0$ as $s\to 0^+$, and that $G(s) \cong \sqrt s$ as $s\to +\infty$ (such behavior at infinity will be important in the sequel). 
We now integrate the inequality \eqref{ine4} between $s$ and $t$, obtaining
\[
\psi(\lambda,s) \le \psi(\lambda,t)   +\lambda \int_{s}^t  G\left(\frac{1}{\lambda^2 \tau}\right) d\tau.
\]
Using \eqref{Jensen},  we infer
\begin{equation*}
P_s(\lambda f) \le  \lambda \psi(\lambda,t)  + \lambda^2 \int_s^{t}  G\left(\frac{1}{\lambda^2 \tau}\right) d\tau.
\end{equation*}
Letting $s\to 0^+$ we conclude
\begin{equation}\label{ine5}
\lambda f \le  \lambda \psi(\lambda,t)  + \lambda^2 \int_0^{t}  G\left(\frac{1}{\lambda^2 \tau}\right) d\tau.
\end{equation}

At this point we   let $B = B(x,r) = \{x\in \M\mid d(y,x)<r\}$, and consider the function $f(y) = - d(y,x)$. Since we clearly have
 \[
 e^{\lambda f} \le e^{-\lambda r} \mathbf 1_{B^c} + \mathbf 1_B,
 \]
 it follows that for every $t>0$ one has
 \[
 e^{\lambda \psi(\lambda,t)(x)}  = P_t(e^{\lambda f})(x) \le e^{-\lambda r} + P_t(\mathbf 1_B)(x).
 \]
 This gives the lower bound
 \[
  P_t(\mathbf 1_B)(x) \ge  e^{\lambda \psi(\lambda,t)(x)} - e^{-\lambda r} .
  \]
To estimate the first term in the right-hand side of the latter inequality, we use \eqref{ine5}  which gives
\[
1 = e^{\lambda f(x)} \le  e^{\lambda \psi(\lambda,t)(x)} e^{\Phi(\lambda,t)},
\]
where we have set
\[
\Phi(\lambda,t) =  \lambda^2 \int_0^{t}  G\left(\frac{1}{\lambda^2 \tau}\right) d\tau.
\]
This gives
\[
 P_{t}(\mathbf 1_B)(x) \ge  e^{-\Phi(\lambda,t)} - e^{-\lambda r}.
 \]
To make use of this estimate, we now choose $\lambda = \frac{1}{r}$, $t = Ar^2$, obtaining
\[
 P_{Ar^2}(\mathbf 1_B)(x) \ge  e^{-\Phi(\frac{1}{r},Ar^2)} - e^{-1}.
\]
We want  to show that we can choose $A>0$ sufficiently small, depending only on $n$, and a $K>0$ (we can in fact take $K = (1-e^{-1})/2$, such that
\begin{equation}\label{ine6}
e^{-\Phi(\frac{1}{r},Ar^2)} - e^{-1} \ge K, \ \ \ \ \text{for every}\ x\in \M, r>0. 
\end{equation}
Consider the function
\[
\Phi(\frac{1}{r},Ar^2) = \frac{1}{r^2} \int_0^{Ar^2}  G\left(\frac{r^2}{\tau}\right) d\tau = \int_{A^{-1}}^\infty \frac{G(t)}{t^2} dt.
\]
As we have noted above a direct examination shows that $G(t) = O(t^{1/2})$ as $t\to +\infty$, and therefore that $\frac{G(t)}{t^2}\in L^1(1,\infty)$. It is then clear that 
\[
\underset{A\to 0^+} {\lim} \int_{A^{-1}}^\infty \frac{G(t)}{t^2} dt = 0,
\]
and therefore there exists $A>0$ sufficiently small such that \eqref{ine6} hold with, say, $K= (1-e^{-1})/2>0$.

\end{proof}

\section{Volume growth and doubling condition}\label{S:doubling}

In this section we prove the following result. 

\begin{theorem}\label{T:doubling}
Let $\M$ be a complete $n$-dimensional Riemannian manifold such that \emph{Ric} $\ge 0$.
There exists a constant $C(n)>0$ such that for every $x\in \M$ and every $r>0$ one has
\[
\emph{Vol}_n(B(x,2r))\le C(n) \emph{Vol}_n(B(x,r)).
\]
\end{theorem}

Of course, Theorem \ref{T:doubling} is contained in the well-known Bishop-Gromov comparison theorem. The point of our proof is that it is purely analytical and does not rely on the analysis of Jacobi fields. Thus, it is possible that our approach might be generalized to the sub-Riemannian setting of \cite{baudoin-garofalo}.  
We will need the following basic result which is a straightforward consequence of the Li Yau inequality (\ref{LiYau}), see \cite{LY}.

\begin{theorem}\label{T:harnack}
Let $\M$ be a $n$-dimensional complete Riemannian manifold such that
\emph{Ric} $\ge 0$. Let $f\in
C^\infty(\bM)$ be such that $0\le f\le M$, and consider $u(x,t) =
P_t f(x)$. For every $(x,s), (y,t)\in \bM\times (0,\infty)$ with
$s<t$ one has 
\begin{equation}\label{beauty}
u(x,s) \le u(y,t) \left(\frac{t}{s}\right)^{\frac{n}{2}} \exp\left(\frac{d(x,y)^2}{4(t-s)} \right).
\end{equation}
\end{theorem}

\begin{corollary}\label{C:harnackheat}
Let $p(x,y,t)$ be the heat kernel on $\bM$. For every $x,y, z\in
\bM$ and every $0<s<t<\infty$ one has
\[
p(x,y,s) \le p(x,z,t) \left(\frac{t}{s}\right)^{\frac{n}{2}}
\exp\left( \frac{d(y,z)^2}{4(t-s)} \right).
\]
\end{corollary}

We now turn to the proof of Theorem \ref{T:doubling}.

\begin{proof}
From the semigroup property and the symmetry of the heat kernel we
have for any $y\in \bM$ and $t>0$
\[ p(y,y,2t) = \int_\bM  p(y,z,t)^2 d\mu(z).
\]

Consider now a function $h\in C^\infty_0(\bM)$ such that $0\le h\le
1$, $h\equiv 1$ on $B(x,\sqrt{t}/2)$ and $h\equiv 0$ outside
$B(x,\sqrt t)$. We thus have
\begin{align*}
P_t h(y) & = \int_\bM p(y,z,t) h(z) d\mu(z) \le \left(\int_{B(x,\sqrt t)}
p(y,z,t)^2 d\mu(z)\right)^{\frac{1}{2}} \left(\int_\bM h(z)^2
d\mu(z)\right)^{\frac{1}{2}}
\\
& \le p(y,y,2t)^{\frac{1}{2}} \mu(B(x,\sqrt t))^{\frac{1}{2}}.
\end{align*}
If we take $y=x$, and $t =r^2$, we obtain
\begin{equation}\label{ine7}
P_{r^2} \left(\mathbf 1_{B(x,r)}\right)(x)^2 \le P_{r^2} h(x)^2 \leq
p(x,x,2r^2)\ \mu(B(x,r)).
\end{equation}
At this point we use the crucial inequality \eqref{doubling}, which gives for some $0<A<1$, depending on $n$,
\[
 P_{Ar^2}(\mathbf 1_{B(x,r)})(x) \ge K, \ \ \ \ \ x\in \M, r>0.
\]
Combining the latter inequality with Theorem \ref{T:harnack} and with \eqref{ine7}, we obtain the following 
on-diagonal lower bound
\begin{equation}\label{odlb}
p(x,x,2r^2) \ge \frac{K^*}{\mu(B(x,r))},\ \ \ \ \ x\in \bM,\ r>0.
\end{equation}
Applying Corollary \ref{C:harnackheat} to $(y,t)\to p(x,y,t)$
for every $y\in B(x,\sqrt t)$ we find
\[
p(x,x,t) \le  C(n) p(x,y,2t).
\]
Integration over $B(x,\sqrt t)$ gives
\[
p(x,x,t)\mu(B(x,\sqrt t)) \le C(n) \int_{B(x,\sqrt
t)}p(x,y,2t)d\mu(y) \le C(n),
\]
where we have used $P_t1\le 1$. Letting $t = r^2$, we obtain from this the on-diagonal upper
bound \begin{equation}\label{odub} p(x,x,r^2) \le
\frac{C(n)}{\mu(B(x,r))}.
\end{equation}
Combining \eqref{odlb} with \eqref{odub} we finally obtain
\[
\mu(B(x,2r)) \le \frac{C}{p(x,x,4r^2)} \le \frac{C^*}{p(x,x,2r^2)}
\le C^{**} \mu(B(x,r)),
\]
where we have used once more Corollary \ref{C:harnackheat}
which gives
\[
\frac{p(x,x,2r^2)}{p(x,x,4r^2)}\le C.
\]
This completes the proof.

\end{proof}

\end{document}